\documentclass[12pt,reqno]{amsart}

\usepackage{geometry}
\usepackage{hhline}

\usepackage{mathtools}

\usepackage{mathdots}

\usepackage{color}

\usepackage{pdfsync}

\usepackage{enumitem}

\usepackage{wasysym}

\usepackage{amssymb}

\usepackage{mathrsfs}

\DeclareMathAlphabet{\mathpzc}{OT1}{pzc}{m}{it}

\usepackage[all]{xy}

\usepackage{tikz}
\usetikzlibrary{arrows,matrix,decorations.pathmorphing,decorations.pathreplacing,positioning,shapes.geometric,shapes.misc,decorations.markings,decorations.fractals,calc,patterns}

\usepackage{graphicx}
\usepackage{float}
\usepackage[bottom]{footmisc}
\usepackage{moreenum}
\usepackage{makecell}
\entrymodifiers={+!!<0pt,\fontdimen22\textfont2>}

\usepackage{scalerel,stackengine}
\stackMath
\newcommand\newcheck[1]{%
\savestack{\tmpbox}{\stretchto{%
  \scaleto{%
    \scalerel*[\widthof{\ensuremath{#1}}]{\kern-.6pt\bigwedge\kern-.6pt}%
    {\rule[-\textheight/2]{1ex}{\textheight}}%WIDTH-LIMITED BIG WEDGE
  }{\textheight}% 
}{0.5ex}}%
\stackon[1pt]{#1}{\scalebox{-1}{\tmpbox}}%
}
\stackMath
\newcommand\newhat[1]{%
\savestack{\tmpbox}{\stretchto{%
  \scaleto{%
    \scalerel*[\widthof{\ensuremath{#1}}]{\kern-.6pt\bigwedge\kern-.6pt}%
    {\rule[-\textheight/2]{1ex}{\textheight}}%WIDTH-LIMITED BIG WEDGE
  }{\textheight}% 
}{0.5ex}}%
\stackon[1pt]{#1}{\scalebox{1}{\tmpbox}}%
}

\def\cA{\mathscr{A}}
\def\cB{\mathscr{B}}
\def\cC{\mathscr{C}}

\def\cF{\mathscr{F}}

\def\cT{\mathscr{T}}

\mathchardef\mhyphen="2D

\def\adots{\mathinner{\mkern1mu\raise1.0pt\vbox{\kern7.0pt\hbox{.}}\mkern2mu\raise5.0pt\hbox{.}\mkern2mu\raise9.0pt\hbox{.}\mkern1mu}}

\def\Coker{\operatorname{Coker}}

\def\dddots{\mathinner{\mkern1mu\raise10.0pt\vbox{\kern7.0pt\hbox{.}}\mkern2mu\raise5.3pt\hbox{.}\mkern2mu\raise1.0pt\hbox{.}\mkern1mu}}
\def\dddotssmall{\mathinner{\mkern1mu\raise7.0pt\vbox{\kern7.0pt\hbox{.}}\mkern-1mu\raise4pt\hbox{.}\mkern-1mu\raise1.0pt\hbox{.}\mkern1mu}}

\def\dual{\operatorname{D}}

\def\id{\operatorname{id}}

\def\K{\operatorname{K}}
\def\K0{\operatorname{K}_0}

\def\Ker{\operatorname{Ker}}

\def\mod{\operatorname{mod}}

\def\PSL2{\operatorname{PSL}_2}

\def\SL2{\operatorname{SL}_2}

\def\TP{\mathsf{TorsPair}}

\numberwithin{equation}{section}

\renewcommand{\labelenumi}{(\roman{enumi})}

\newtheorem{Lemma}{Lemma}[section]
\theoremstyle{definition}
\newtheorem{Definition}[Lemma]{Definition}

\newtheorem{Remark}[Lemma]{Remark}

\newtheorem{ThmIntro}{Theorem}

\newtheorem*{bfhpg*}{}

  {\begin{list}{}{%
    \settowidth{\labelwidth}{\textbf{#1:}}%
    \setlength{\leftmargin}{\labelwidth}\addtolength{\leftmargin}{\labelsep}}}%
  {\end{list}}

\makeatletter
\@namedef{subjclassname@2020}{%
  \textup{2020} Mathematics Subject Classification}
\makeatother

\begin{document}

\setlength{\parindent}{0pt}
\setlength{\parskip}{7pt}

\title[Proper abelian subcategories and tilting]{Proper abelian subcategories of triangulated categories and their tilting theory}

\author{Peter J\o rgensen}

\address{Department of Mathematics, Aarhus University, Ny Munkegade 118, 8000 Aarhus C, Denmark}
\email{peter.jorgensen@math.au.dk}

\urladdr{https://sites.google.com/view/peterjorgensen}

%\thanks{Date: \today. A thank you would go here}

\keywords{Happel--Reiten--Smal\o\ tilting, heart, $t$-structure, torsion pair}

\subjclass[2020]{18E10, 18E40, 18G80}

%05B45: Combinatorial aspects of tessellation and tiling problems
%05E10: Combinatorial aspects of representation theory
%05E16: Combinatorial aspects of groups and algebras
%05E40: Combinatorial aspects of commutative algebra
%05E45: Combinatorial aspects of simplicial complexes
%05E99: Algebraic combinatorics / None of the above, but in this section
%13F60: Cluster algebras
%16E10: Homological dimension in associative algebras
%16E45: Differential graded algebras and applications (associative algebraic aspects)
%16G10: Representations of associative Artinian rings 
%16G60: Representation type (finite, tame, wild, etc.) of associative algebras
%16G70: Auslander-Reiten sequences (almost split sequences) and Auslander-Reiten quivers
%16S90: Torsion theories; radicals on module categories (associative algebraic aspects) {For radicals of rings, see 16Nxx}
%18A20: Epimorphisms, monomorphisms, special classes of morphisms, null morphisms
%18E10: Abelian categories, Grothendieck categories
%18E35: Localization of categories, calculus of fractions
%18E40: Torsion theories, radicals
%18G05: Projectives and injectives (category-theoretic aspects) [See also 13C10, 13C11, 16D40, 16D50]
%18G35: Chain complexes (category-theoretic aspects), dg categories [See also 14F07, 18G80, 55U15]
%18G80: Derived categories, triangulated categories
%18G99: Homological algebra in category theory, derived categories and functors / None of the above, but in this section 
%51M20: Polyhedra and polytopes; regular figures, division of spaces
%55P62: Rational homotopy theory

\begin{abstract} 

In the theory of triangulated categories, we propose to replace {\em hearts of $t$-structures} by {\em proper abelian subcategories}, which may be plentiful even when hearts are not.  For instance, this happens in negative cluster categories.  

\medskip
\noindent
In support of our proposal, we show that proper abelian subcategories with a few vanishing negative self extensions permit a tilting theory which is a direct generalisation of Happel--Reiten--Smal\o\ tilting of hearts.

\end{abstract}

\maketitle

\setcounter{section}{-1}
\section{Introduction}
\label{sec:introduction}

The purpose of this paper is to propose a change of perspective from {\em hearts of $t$-structures} to {\em proper abelian subcategories} of triangulated categories.

\begin{Definition}
[{see \cite[def.\ 1.2]{JorSMS}}]
\label{def:proper_abelian_subcategory}
A {\em proper abelian subcategory} $\cA$ of a triangulated category $\cC$ is an additive subcategory which is abelian and satisfies that $0 \xrightarrow{} a' \xrightarrow{ \alpha' } a \xrightarrow{ \alpha } a'' \xrightarrow{} 0$ is a short exact sequence in $\cA$ if and only if $a'$, $a$, $a''$ are in $\cA$ and there is a triangle $a' \xrightarrow{ \alpha' } a \xrightarrow{ \alpha } a'' \xrightarrow{} \Sigma a'$ in $\cC$.  
\end{Definition}

Hearts of $t$-structures are proper abelian subcategories, but there are many others.  Indeed, the heart of a $t$-structure has all negative self extensions equal to zero, but this property fails for large classes of proper abelian subcategories arising in practice.  For instance, \cite[Introduction]{JorSMS} provides an explanation of this in the setting of negative cluster categories as developed in \cite{CS2}, \cite{CS}, \cite{CSP1}, \cite{CSPP}, \cite{IJ}, \cite{Jin}.

In support of the proposal to supplant hearts of $t$-structures, we will show that proper abelian subcategories with a few vanishing negative self extensions permit a tilting theory which is a direct generalisation of Happel--Reiten--Smal\o\ tilting of hearts, see \cite[thm.\ I.3.1]{BR}, \cite[prop.\ 2.1]{HRS}, and \cite[prop.\ 2.1]{W}.  To state our results, we need the following vanishing condition.

\begin{Definition}
\label{def:E}
Let $n \geqslant 1$ be an integer or $\infty$.  A full subcategory $\cA$ of a triangulated category $\cC$ {\em satisfies condition $E_n$} if $\cC( \cA,\Sigma^i\cA ) = 0$ for $-n-1 < i < 0$.
\end{Definition}

Observe that if $n \geqslant 2$ then condition $E_n$ implies condition $E_{ n-1 }$.  We also need the following generalisation of the class of right tilts of a fixed heart.

\begin{Definition}
\label{def:intermediate}
Let $\cA$ be a proper abelian subcategory of a triangulated category $\cC$ and let $n \geqslant 1$ be an integer or $\infty$.  The {\em $E_n$-interval right of $\cA$ in $\cC$} is
\[
  [ \cA,\Sigma \cA ]_n
  =
  \left\{\!\!
      \begin{array}{ll}
        \mbox{ \rm proper abelian subcategories $\cB \subseteq \cC$ satisfying $E_n$ } \\[1.5mm]
        \mbox{ \rm such that $\cA \subseteq \cB * \Sigma^{ -1 }\cB$ and $\cB \subseteq \Sigma \cA * \cA$ }\\
      \end{array}
    \!\right\}.
\]
For $\cB', \cB \in [ \cA,\Sigma \cA ]_n$ we write
\[
  \cB' \preccurlyeq \cB
    \;\; \mbox{ if } \;\;
  \cB' \subseteq \cB * \Sigma^{ -1 }\cB
    \mbox{ and }
  \cB \subseteq \Sigma \cB' * \cB'.
\]
\end{Definition}

Under mild conditions, $\preccurlyeq$ is a partial order on $[ \cA,\Sigma \cA ]_n$ which generalises the partial order on right tilts of a fixed heart given by inclusion of the corresponding coailes.

\begin{ThmIntro}
\label{thm:partial_order}
Let $\cA$ be a proper abelian subcategory of a Krull--Schmidt triangulated category $\cC$ and let $n \geqslant 3$ be an integer or $\infty$.
\begin{enumerate}
\setlength\itemsep{4pt}

  \item  $\preccurlyeq$ is a partial order on $[ \cA,\Sigma \cA ]_n$.

  \item  If $\cA$ satisfies condition $E_n$ then $\cA$ and $\Sigma \cA$ are the least and greatest elements for $\preccurlyeq$.

\end{enumerate}
\end{ThmIntro}

Our main result is the following direct generalisation of Happel--Reiten--Smal\o\ tilting.

\begin{ThmIntro}
\label{thm:proper_tilting13}
Let $\cA$ be a proper abelian subcategory of a Krull--Schmidt triangulated category $\cC$, let $n \geqslant 4$ be an integer or $\infty$, and assume that $\cA$ satisfies condition $E_n$.  Then there are order preserving, mutually inverse bijections
\[
  \xymatrix {
  \TP( \cA )
    \ar[rr]<1ex>^-{ \Phi( - ) } &&
  [ \cA,\Sigma \cA ]_{ n-1 }
    \ar[ll]<1ex>^-{ \Psi( - ) }
            }
\]
given by
\[
  \Phi( \cT,\cF ) = \Sigma \cF * \cT
  \;\;,\;\;
  \Psi( \cB ) = ( \cA \cap \cB,\cA \cap \Sigma^{ -1 }\cB ).
\]
\end{ThmIntro}

Here $\TP( \cA )$ is the class of torsion pairs in $\cA$ and $*$ is defined as usual, see Equation \eqref{equ:star}.  The partial orders preserved by the bijections are $\sqsubseteq$ on $\TP( \cA )$, defined by inclusion of torsion free classes, and $\preccurlyeq$ on $[ \cA,\Sigma \cA ]_{ n-1 }$, see Definition \ref{def:intermediate}.

In Theorem \ref{thm:proper_tilting13}, the tilting map $\Phi$ starts with a torsion pair $( \cT,\cF )$ in a proper abelian subcategory $\cA$ satisfying condition $E_n$, but produces a proper abelian subcategory $\Phi( \cT,\cF )$ satisfying condition $E_{ n-1 }$.  The following result shows that this can sometimes be improved, hence permitting the tilting process to be iterated indefinitely.

\begin{ThmIntro}
\label{thm:proper_tilting16}
In Theorem \ref{thm:proper_tilting13}, make either of the following additional assumptions.
\begin{enumerate}
\setlength\itemsep{4pt}

  \item  $n = \infty$.

  \item  $n \geqslant 4$ is an integer and $\cC$ is $k$-linear $( -n-1 )$-Calabi--Yau over a field $k$.
  
\end{enumerate}
Then there are mutually inverse bijections
\[
  \xymatrix {
  \TP( \cA )
    \ar[rr]<1ex>^-{ \Phi( - ) } &&
  [ \cA,\Sigma \cA ]_n
    \ar[ll]<1ex>^-{ \Psi( - ) }
            }
\]
given by the formulae in Theorem \ref{thm:proper_tilting13}.
\end{ThmIntro}

In particular, if $Q$ is a finite acyclic quiver, $k$ a field, then part (ii) of Theorem \ref{thm:proper_tilting16} applies to the ``generic'' proper abelian subcategory $\cA = \mod kQ$ of the negative cluster category $\cC_{ -n-1 }( kQ )$, see \cite[sec.\ 1.2]{CSPP}.  The proofs of Theorems \ref{thm:proper_tilting13} and \ref{thm:proper_tilting16} are facilitated by the following simple criterion for proper abelian subcategories, which has the classic \cite[prop.\ 1.2.4]{BBD} from the theory of $t$-structures as a special case.

\begin{ThmIntro}
\label{thm:SMS7}
Let $\cA$ be an additive, extension closed subcategory of a triangulated category $\cC$, and assume $\cA$ satisfies condition $E_2$.  Then 
\[
  \mbox{ $\cA$ is a proper abelian subcategory of $\cC$ $\Leftrightarrow$ $\cA * \Sigma \cA \subseteq \Sigma \cA * \cA$. }
\]
\end{ThmIntro}

Finally, note that related work has recently appeared in \cite{Klapproth} and \cite{L}.  The notion of {\em proper $n$-exact subcategories} was in effect considered in \cite[thm.\ 1]{Klapproth}, and the notion of {\em distinguished abelian subcategories} was introduced in \cite[def.\ 1.1]{L}.

\section{Setup}
\label{sec:setup}

The following applies throughout the paper.   
\begin{itemize}
\setlength\itemsep{4pt}

  \item  Consider an additive category.  An {\em additive subcategory} is a full subcategory closed under isomorphisms, direct sums, and direct summands.  

  \item  The suspension functor of a triangulated category is denoted $\Sigma$. 
  
  \item  A {\em short triangle} in a triangulated category is a diagram $c' \xrightarrow{ \gamma' } c \xrightarrow{ \gamma } c''$ for which there is a triangle $c' \xrightarrow{ \gamma' } c \xrightarrow{ \gamma } c'' \xrightarrow{ \gamma'' } \Sigma c'$.

  \item  Let $\cA$, $\cB$ be full subcategories of a triangulated category $\cC$.  Then
\begin{equation}
\label{equ:star}
  \cA * \cB = \{ e \in \cC \mid \mbox{there is a short triangle $a \xrightarrow{} e \xrightarrow{} b$ with $a \in \cA$ and $b \in \cB$} \}
\end{equation}
is a full subcategory of $\cC$ closed under isomorphisms.

  \item  If $\cA$ and $\cB$ are closed under direct sums, then so is $\cA * \cB$.
  
  \item  If $\cC$ is Krull--Schmidt and $\cA$ and $\cB$ are additive subcategories satisfying $\cC( \cA,\cB ) = 0$, then $\cA * \cB$ is an additive subcategory of $\cC$ by \cite[prop.\ 2.1(1)]{IY}.  Note that the proof of \cite[prop.\ 2.1(1)]{IY} only requires the Krull--Schmidt property of $\cC$, not the stronger blanket assumptions of \cite{IY}.
  
  \item  We say that $\cA$ is {\em extension closed in $\cC$} if $\cA * \cA \subseteq \cA$.  

\end{itemize}

%\begin{Remark}
%Although $[ \cA,\Sigma \cA ]_n$ depends on $\cC$ as well as $\cA$, we have omitted $\cC$ from the notation to avoid making it too heavy.
%\end{Remark}

%\begin{Remark}
%Definition \ref{def:E} could have been formulated in the following two parts:
%\begin{itemize}
%\setlength\itemsep{4pt}
%
%  \item  If $n \geqslant 1$ is an integer, then $\cA$ satisfies condition $E_n$ if $\cC( \cA,\Sigma^i \cA ) = 0$ for $i \in \{ -1, \ldots, -n \}$.
%
%  \item  $\cA$ satisfies condition $E_{ \infty }$ if $\cC( \cA,\Sigma^i \cA ) = 0$ for $i \leqslant -1$.
%
%\end{itemize}
%We chose the shorter Definition \ref{def:E} because it makes some of the proofs shorter too.
%\end{Remark}

\section{Lemmas on additive subcategories}
\label{sec:lemmas1}

\begin{Lemma}
\label{lem:pre_proper_tilting8a}
Let $\cA$ be an additive subcategory of a triangulated category $\cC$ and assume $\cA$ satisfies condition $E_1$.  Then $\cA * \Sigma^2 \cA \subseteq \Sigma^2 \cA * \cA$.  
\end{Lemma}

\begin{proof}
Given $c \in \cA * \Sigma^2 \cA$, there is a triangle $a' \xrightarrow{} c \xrightarrow{} \Sigma^2 a'' \xrightarrow{ \delta } \Sigma a'$ with $a', a'' \in \cA$.  Here $\delta = 0$ since $\cA$ satisfies condition $E_1$, so the triangle is split whence $c \cong a' \oplus \Sigma^2 a''$.  Hence there is also a split triangle $\Sigma^2 a'' \xrightarrow{} c \xrightarrow{} a' \xrightarrow{ 0 } \Sigma\Sigma^2 a''$ which proves $c \in \Sigma^2 \cA * \cA$.  
\end{proof}

\begin{Lemma}
\label{lem:SMS6}
Let $\cA$ be an additive subcategory of a triangulated category $\cC$ and assume $\cA$ satisfies condition $E_2$.
\begin{enumerate}
\setlength\itemsep{4pt}

  \item  $( \Sigma\cA * \cA ) \cap ( \cA * \Sigma^{ -1 } \cA ) = \cA$.

  \item  If $\cA * \Sigma \cA \subseteq \Sigma \cA * \cA$ then $( \Sigma^{ -1 }\cA * \cA ) \cap ( \cA * \Sigma \cA ) = \cA$.

\end{enumerate}
\end{Lemma}

\begin{proof}
(i):  The inclusion $\supseteq$ is clear.  To prove $\subseteq$, let $c \in ( \Sigma\cA * \cA ) \cap ( \cA * \Sigma^{ -1 }\cA )$ be given.  Since $c \in \cA * \Sigma^{ -1 }\cA$, there is a short triangle $a' \xrightarrow{} c \xrightarrow{ \gamma } \Sigma^{ -1 }a''$ with $a', a'' \in \cA$.  Since $c \in \Sigma\cA * \cA$, condition $E_2$ implies $\gamma = 0$ whence $c$ is a direct summand of $a'$ so $c \in \cA$ as desired.

(ii):  The inclusion $\supseteq$ is clear, and $\subseteq$ can be proved as follows:
\[
  ( \Sigma^{ -1 }\cA * \cA ) \cap ( \cA * \Sigma \cA )
  \subseteq 
  ( \cA * \Sigma^{ -1 }\cA ) \cap ( \Sigma \cA * \cA )
  =
  \cA,
\]
where we used $\cA * \Sigma \cA \subseteq \Sigma \cA * \cA$, its first negative suspension, and part (i).
\end{proof}

\begin{Lemma}
\label{lem:pre_proper_tilting8b}
Let $\cA$ be an additive subcategory of a Krull--Schmidt triangulated category $\cC$ and assume $\cA$ satisfies condition $E_2$.  Then $( \cA * \Sigma \cA )\cap ( \Sigma^2 \cA * \Sigma \cA * \cA * \Sigma^{ -1 }\cA ) \subseteq \Sigma \cA * \cA$.
\end{Lemma}

\begin{proof}
First observe that condition $E_2$ and \cite[prop.\ 2.1(1)]{IY} imply that $\Sigma \cA * \cA$ and $\Sigma \cA * \cA * \Sigma^{ -1 }\cA$ are additive subcategories of $\cC$.

Let $c \in ( \cA * \Sigma \cA )\cap ( \Sigma^2 \cA * \Sigma \cA * \cA * \Sigma^{ -1 }\cA )$ be given.  Since $c \in \Sigma^2 \cA * \Sigma \cA * \cA * \Sigma^{ -1 }\cA$, there is a short triangle $\Sigma^2 a' \xrightarrow{ \sigma } c \xrightarrow{} x$ with $a' \in \cA$ and $x \in \Sigma \cA * \cA * \Sigma^{ -1 }\cA$.  Since $c \in \cA * \Sigma \cA$, condition $E_2$ implies $\sigma = 0$ whence $c$ is a direct summand of $x$ so $c \in \Sigma \cA * \cA * \Sigma^{ -1 }\cA$ since $\Sigma \cA * \cA * \Sigma^{ -1 }\cA$ is an additive subcategory.

Hence there is a short triangle $y \xrightarrow{} c \xrightarrow{ \gamma } \Sigma^{ -1 }a''$ with $y \in \Sigma \cA * \cA$ and $a'' \in \cA$.  Since $c \in \cA * \Sigma \cA$, condition $E_2$ implies $\gamma = 0$ whence $c$ is a direct summand of $y$ so $c \in \Sigma \cA * \cA$ as desired since $\Sigma \cA * \cA$ is an additive subcategory.
\end{proof}

\begin{Lemma}
\label{lem:proper_tilting2}
Let $\cA$ be an additive subcategory of a Krull--Schmidt triangulated category $\cC$ and assume $\cA$ satisfies condition $E_3$.  Then $( \Sigma \cA * \cA * \Sigma^{ -1 } \cA ) \cap ( \cA * \Sigma^{ -1 } \cA * \Sigma^{ -2 } \cA ) = \cA * \Sigma^{ -1 } \cA$.  
\end{Lemma}

\begin{proof}
The inclusion $\supseteq$ is clear.  To prove $\subseteq$, let $c \in ( \Sigma \cA * \cA * \Sigma^{ -1 } \cA ) \cap ( \cA * \Sigma^{ -1 } \cA * \Sigma^{ -2 } \cA )$ be given.  Since $c \in \Sigma \cA * \cA * \Sigma^{ -1 } \cA$, there is a short triangle $\Sigma a \xrightarrow{ \sigma } c \xrightarrow{} x$ with $a \in \cA$ and $x \in \cA * \Sigma^{ -1 } \cA$.  Since $c \in \cA * \Sigma^{ -1 } \cA * \Sigma^{ -2 }\cA$, condition $E_3$ implies $\sigma = 0$.  Hence $c$ is a direct summand of $x$.  However, condition $E_3$ and \cite[prop.\ 2.1(1)]{IY} imply that $\cA * \Sigma^{ -1 } \cA$ is an additive subcategory (actually $E_1$ suffices), so since $x$ is in $\cA * \Sigma^{ -1 } \cA$, so is $c$.  
\end{proof}

\section{Proof of Theorem \ref{thm:partial_order}}
\label{sec:partial_order}

\begin{proof}
[Proof of Theorem \ref{thm:partial_order}]
(i):  It is clear that $\preccurlyeq$ is reflexive.

To prove that $\preccurlyeq$ is transitive, let $\cB'' \preccurlyeq \cB' \preccurlyeq \cB$ in $[ \cA,\Sigma \cA ]_n$ be given.  On the one hand we have
\[
  \cB''
  \subseteq
  \cB' * \Sigma^{ -1 }\cB'
  \subseteq
  \cB * \Sigma^{ -1 }\cB * \Sigma^{ -1 }\cB * \Sigma^{ -2 }\cB
  = \cB * \Sigma^{ -1 }\cB * \Sigma^{ -2 }\cB.
\]
On the other hand we have
\[
  \cB''
  \subseteq
  \Sigma \cA * \cA
  \subseteq
  \Sigma \cB * \cB * \cB * \Sigma^{ -1 }\cB
  \subseteq
  \Sigma \cB * \cB * \Sigma^{ -1 }\cB.
\]
The computations use that $\Sigma^{ -1 } \cB$ and $\cB$ are closed under extensions.  Lemma \ref{lem:proper_tilting2} now shows 
\[
%\label{equ:partial_order1}
  \cB'' \subseteq \cB * \Sigma^{ -1 }\cB.
\]
Similar computations show
%Similarly, on the one hand we have
%\[
%  \cB
%  \subseteq
%  \Sigma \cB' * \cB'
%  \subseteq
%  \Sigma^2 \cB'' * \Sigma \cB'' * \Sigma \cB'' * \cB''
%  = \Sigma^2 \cB'' * \Sigma \cB'' * \cB''.
%\]
%On the other hand we have
%\[
%  \cB
%  \subseteq
%  \Sigma \cA * \cA
%  \subseteq
%  \Sigma \cB'' * \cB'' * \cB'' * \Sigma^{ -1 }\cB''
%  \subseteq
%  \Sigma \cB'' * \cB'' * \Sigma^{ -1 }\cB''.
%\]
%The two computations use that $\Sigma \cB''$ and $\cB''$ are closed under extensions.  The suspension of Lemma \ref{lem:proper_tilting2} now shows 
\[
  \cB \subseteq \Sigma \cB'' * \cB'',
\]
and 
%combined with Equation \eqref{equ:partial_order1} this 
the last two equations
show $\cB'' \preccurlyeq \cB$ as desired.

To prove that $\preccurlyeq$ is antisymmetric, let $\cB', \cB$ in $[ \cA,\Sigma \cA ]_n$ satisfy $\cB' \preccurlyeq \cB$ and $\cB \preccurlyeq \cB'$.  Then $\cB' \preccurlyeq \cB$ gives $\cB' \subseteq \cB * \Sigma^{ -1 }\cB$ and $\cB \subseteq \Sigma \cB' * \cB'$ while $\cB \preccurlyeq \cB'$ gives $\cB' \subseteq \Sigma \cB * \cB$ and $\cB \subseteq \cB' * \Sigma^{ -1 }\cB'$.  Lemma \ref{lem:SMS6}(i) now gives $\cB' \subseteq \cB$ and $\cB \subseteq \cB'$ whence $\cB' = \cB$.

(ii):  Assume that $\cA$ satisfies condition $E_n$.  Then it is immediate from Definition \ref{def:intermediate} that $\cA, \Sigma \cA \in [ \cA,\Sigma \cA ]_n$ and that each $\cB \in [ \cA,\Sigma \cA ]_n$ satisfies $\cA \preccurlyeq \cB \preccurlyeq \Sigma \cA$.
\end{proof}

\section{Lemmas on proper abelian subcategories}
\label{sec:lemmas2}

The following lemma was proved in \cite[lem.\ 4.2]{JorSMS}.  Note that the proof only requires the assumptions made here, not the stronger blanket assumptions of \cite{JorSMS}.

\begin{Lemma}
\label{lem:SMS23}
Let $\cA$ be a proper abelian subcategory of a triangulated category.  If $a_1 \xrightarrow{ \alpha_1 } a_0 \xrightarrow{} c$ is a short triangle with $a_1, a_0 \in \cA$, then there is a short triangle $\Sigma \Ker \alpha_1 \xrightarrow{} c \xrightarrow{} \Coker \alpha_1$ where $\Ker \alpha_1$ and $\Coker \alpha_1$ are the kernel and cokernel of $\alpha_1$ in $\cA$.
\end{Lemma}

%\begin{proof}
%The morphism $\alpha_1$ can be factored as $a_1 \xrightarrow{ q } x \xrightarrow{ i } a_0$ where $x$ is the image of $\alpha_1$ in $\cA$ while $q$ and $i$ are an epimorphism and a monomorphism in $\cA$.  There are short exact sequences $a' \xrightarrow{} a_1 \xrightarrow{ q } x$ and $x \xrightarrow{ i } a_0 \xrightarrow{} a''$ in $\cA$ with $a' = \Ker \alpha_1$ and $a'' = \Coker \alpha_1$.  The sequences are short triangles in $\cC$, and they can be combined by the octahedral axiom to the following commutative diagram where each row and column is a short triangle; see \cite[prop.\ 1.4.6]{N}.
%\begin{equation}
%\label{equ:SMS23_1}
%\vcenter{
%  \xymatrix @+0.5pc {
%    a' \ar[r] \ar[d] & a_1 \ar^{ q }[r] \ar@{=}[d] & x \ar^{ i }[d] \\
%    \Sigma^{ -1 }c \ar[r] \ar[d] & a_1 \ar_{ \alpha_1 }[r] \ar[d] & a_0 \ar[d] \\
%    \Sigma^{ -1 }a'' \ar[r] & 0 \ar[r] & a'' \\
%                    }
%        }
%\end{equation}
%The suspension of the first column provides the desired short triangle.
%\end{proof}

\begin{Lemma}
\label{lem:criterion_for_torsion_pair}
Let $\cA$ be a proper abelian subcategory of a triangulated category $\cC$ and let $\cT, \cF$ be full subcategories of $\cA$ closed under isomorphism.  Then $( \cT,\cF )$ is a torsion pair in $\cA$ if and only if $\cC( \cT,\cF ) = 0$ and $\cA = \cT * \cF$.
\end{Lemma}

\begin{proof}
We know by \cite[p.\ 224]{Dickson} that $( \cT,\cF )$ is a torsion pair in $\cA$ if and only if $\cC( \cT,\cF ) = 0$ and each $a \in \cA$ permits a short exact sequence $t \xrightarrow{} a \xrightarrow{} f$ with $t \in \cT$ and $f \in \cF$.  Since $\cA$ is a proper abelian subcategory, the latter condition is equivalent to each $a \in \cA$ permitting a short triangle $t \xrightarrow{} a \xrightarrow{} f$ with $t \in \cT$ and $f \in \cF$, and this is equivalent to $\cA = \cT * \cF$.  
\end{proof}

\begin{Lemma}
\label{lem:proper_tilting567}
Let $\cA$ be a proper abelian subcategory of a triangulated category $\cC$, let $( \cT,\cF )$ be a torsion pair in $\cA$, and consider the full subcategory $\cB = \Sigma \cF * \cT$ of $\cC$.  
\begin{enumerate}
\setlength\itemsep{4pt}

  \item  $\cT$ and $\cF$ are additive, extension closed subcategories of $\cC$.
  
  \item  $\cT * \Sigma \cF \subseteq \Sigma \cF * \cT$.  

  \item  $\cA \subseteq \cB * \Sigma^{ -1 }\cB$ and $\cB \subseteq \Sigma \cA * \cA$.  

  \item  $\cB$ is extension closed in $\cC$.

  \item  If $\cC$ is Krull--Schmidt and $\cA$ satisfies condition $E_1$, then $\cB$ is an additive subcategory of $\cC$. 

  \item  If $\cA$ satisfies condition $E_1$ then $\cA \cap \cB = \cT$ and $\cA \cap \Sigma^{ -1 }\cB = \cF$.

  \item  If $n \geqslant 2$ is an integer or $\infty$ and $\cA$ satisfies condition $E_n$, then $\cB$ satisfies condition $E_{ n-1 }$.   

\end{enumerate}
\end{Lemma}

\begin{proof}
(i) is clear from the definition of proper abelian subcategory and torsion pair.

(ii):  See \cite[lem.\ 4.4(i)]{JorSMS}, the proof of which applies under the assumptions made here.
%Let $c \in \cT * \Sigma \cF$ be given.  There is a short triangle $f \xrightarrow{ \varphi } t \xrightarrow{} c$ with $f \in \cF$, $t \in \cT$, hence a short triangle $\Sigma \Ker \varphi \xrightarrow{} c \xrightarrow{} \Coker \varphi$ by Lemma \ref{lem:SMS23}.  It is hence enough to show $\Ker \varphi \in \cF$ and $\Coker \varphi \in \cT$, and these follow because $\Ker \varphi$ is a subobject of $f \in \cF$ while $\Coker \varphi$ is a quotient object of $t \in \cT$.

(iii):  Using associativity of $*$ which holds by \cite[lem.\ 1.3.10]{BBD}, we have
\[
  \cA \stackrel{ \rm (a) }{ = } \cT * \cF \subseteq \Sigma \cF * \cT * \cF * \Sigma^{ - 1 }\cT = \cB * \Sigma^{ -1 }\cB.
\]
Here (a) is by Lemma \ref{lem:criterion_for_torsion_pair}.  A similar computation proves the second inclusion in (iii).

(iv):  Associativity of $*$ again gives
\[
  \cB * \cB
    =
  \Sigma \cF * \cT * \Sigma \cF * \cT
    \stackrel{ \rm (a) }{ \subseteq }
  \Sigma \cF * \Sigma \cF * \cT * \cT
    \stackrel{ \rm (b) }{ \subseteq }
  \Sigma \cF * \cT
    =
  \cB.
\]
Here (a) and (b) are by parts (ii) and (i).

(v):  Assume that $\cC$ is Krull--Schmidt and that $\cA$ satisfies condition $E_1$.  Then $\cC( \Sigma \cF,\cT ) = 0$ since $\cF, \cT \subseteq \cA$, so $\cB = \Sigma \cF * \cT$ is an additive subcategory of $\cC$ by \cite[prop.\ 2.1(1)]{IY}.

(vi):  To show $\cA \cap \cB = \cT$, note that $\supseteq$ is clear.  For $\subseteq$, let $b \in \cA \cap \cB$ be given.  Since $b \in \cB$, there is a short triangle $\Sigma f \xrightarrow{ \sigma } b \xrightarrow{} t$ with $f \in \cF$, $t \in \cT$.  Since $b \in \cA$ we have $\sigma = 0$ by condition $E_1$ for $\cA$.  Hence $b$ is a direct summand of $t$ whence $b \in \cT$.  A similar argument proves $\cA \cap \Sigma^{ -1 }\cB = \cF$.
%To show $\cA \cap \Sigma^{ -1 }\cB = \cF$, note that $\supseteq$ is clear.  For $\subseteq$, let $\Sigma^{ -1 }b' \in \cA \cap \Sigma^{ -1 }\cB$ be given.  Since $\Sigma^{ -1 }b' \in \Sigma^{ -1 }\cB$, there is a short triangle $f' \xrightarrow{} \Sigma^{ -1 }b' \xrightarrow{ \sigma' } \Sigma^{ -1 }t'$ with $f' \in \cF$, $t' \in \cT$.  Since $\Sigma^{ -1 }b' \in \cA$ we have $\sigma' = 0$ by condition $E_1$ for $\cA$.  Hence $\Sigma^{ -1 }b'$ is a direct summand of $f'$ whence $\Sigma^{ -1 }b' \in \cF$.

(vii):  Assume that $\cA$ satisfies condition $E_n$.  To see that $\cB = \Sigma \cF * \cT$ satisfies condition $E_{ n-1 }$, it is enough to prove the following.
\begingroup
\renewcommand{\labelenumi}{(\alph{enumi})}
\begin{enumerate}
\setlength\itemsep{4pt}

  \item  $\cC( \Sigma \cF,\Sigma^i\Sigma \cF ) = 0$ for $-n < i < 0$.

  \item  $\cC( \Sigma \cF,\Sigma^i\cT ) = 0$ for $-n < i < 0$.
  
  \item  $\cC( \cT,\Sigma^i\Sigma \cF ) = 0$ for $-n < i < 0$.

  \item  $\cC( \cT,\Sigma^i \cT ) = 0$ for $-n < i < 0$.

\end{enumerate}
\endgroup
Since $\cF, \cT \subseteq \cA$, condition $E_n$ for $\cA$ immediately implies parts (a), (b), and (d), and also part (c) for $-n < i < -1$.  Part (c) with $i = -1$ states $\cC( \cT,\cF ) = 0$,  which is true because $( \cT,\cF )$ is a torsion pair in $\cA$.
\end{proof}

\begin{Lemma}
\label{lem:proper_tilting11}
Let $\cA$ be a proper abelian subcategory, $\cB$ an additive subcategory of a triangulated category $\cC$.  Assume that $\cA$ satisfies condition $E_3$, that $\cB$ satisfies condition $E_1$, and that $\cA \subseteq \cB * \Sigma^{ -1 }\cB$ and $\cB \subseteq \Sigma \cA * \cA$.

Then $( \cT,\cF ) := ( \cA \cap \cB,\cA \cap \Sigma^{ -1 }\cB )$ is a torsion pair in $\cA$.  
\end{Lemma}

\begin{proof}
By Lemma \ref{lem:criterion_for_torsion_pair} we must prove $\cC( \cT,\cF ) = 0$ and $\cA = \cT * \cF$, and the former follows from condition $E_1$ for $\cB$.

To prove $\cA = \cT * \cF$, let $a \in \cA$ be given.  Since $\cA \subseteq \cB * \Sigma^{ -1 }\cB$, there is a triangle
\begin{equation}
\label{equ:proper_tilting11a}
  \Sigma^{ -2 }b'' \xrightarrow{} b' \xrightarrow{} a \xrightarrow{} \Sigma^{ -1 }b''
\end{equation}
with $b', b'' \in \cB$, so it is enough to show $b' \in \cT$ and $\Sigma^{ -1 }b'' \in \cF$, hence enough to show $b', \Sigma^{ -1 }b'' \in \cA$.  

To show $b' \in \cA$, note that we have $b' \in \cB \subseteq \Sigma \cA * \cA$, so there is a short triangle $\Sigma a' \xrightarrow{ \sigma' } b' \xrightarrow{} a''$ with $a', a'' \in \cA$.  The triangle \eqref{equ:proper_tilting11a} shows $b' \in \Sigma^{ -2 }\cB * \cA$, so $b' \in \Sigma^{ -1 }\cA * \Sigma^{ -2 }\cA * \cA$ because $\cB \subseteq \Sigma \cA * \cA$.  Condition $E_3$ for $\cA$ hence implies $\sigma' = 0$ so $b'$ is a direct summand of $a'' \in \cA$ whence $b' \in \cA$.  
A similar argument shows $\Sigma^{ -1 }b'' \in \cA$.
%To show $\Sigma^{ -1 }b'' \in \cA$, note that we have $b'' \in \cB \subseteq \Sigma \cA * \cA$, so there is a short triangle $a_0 \xrightarrow{} \Sigma^{ -1 }b'' \xrightarrow{ \sigma'' } \Sigma^{ -1 }a_1$ with $a_0, a_1 \in \cA$.  The triangle \eqref{equ:proper_tilting11a} induces a triangle $a \xrightarrow{} \Sigma^{ -1 }b'' \xrightarrow{} \Sigma b' \xrightarrow{} \Sigma a$ which shows $\Sigma^{ -1 }b'' \in \cA * \Sigma\cB$, so $\Sigma^{ -1 }b'' \in \cA * \Sigma^2\cA * \Sigma \cA$ because $\cB \subseteq \Sigma \cA * \cA$.  Condition $E_3$ for $\cA$ hence implies $\sigma'' = 0$ so $\Sigma^{ -1 }b''$ is a direct summand of $a_0 \in \cA$ whence $\Sigma^{ -1 }b'' \in \cA$.  
\end{proof}

\section{Proof of Theorem \ref{thm:SMS7}}
\label{sec:abelian}

\begin{proof}
[Proof of Theorem \ref{thm:SMS7}]
$\Rightarrow$:  Immediate by Lemma \ref{lem:SMS23}.

$\Leftarrow$:  Assume that $\cA * \Sigma \cA \subseteq \Sigma \cA * \cA$.  Since $\cA \subseteq \cC$ is an additive subcategory satisfying condition $E_2$, we know by \cite[theorem]{Dyer} that $\cA$ is an exact category with conflations given by all short triangles in $\cC$ with terms in $\cA$ (actually $E_1$ suffices).  See \cite[thm.\ 2.5]{JorSMS} for an alternative proof.  By \cite[ex.\ 8.6(i)]{B} it remains to show that each morphism in $\cA$ factors into a deflation followed by an inflation.

Let $a_1 \xrightarrow{ \alpha_1 } a_0$ be a morphism in $\cA$ and complete it to a short triangle $a_1 \xrightarrow{ \alpha_1 } a_0 \xrightarrow{} c$.  Then $c \in \cA * \Sigma \cA \subseteq \Sigma \cA * \cA$ so there is a short triangle $\Sigma a' \xrightarrow{} c \xrightarrow{} a''$ with $a', a'' \in \cA$.  By the octahedral axiom, the short triangles can be combined to the following commutative diagram
%\eqref{equ:SMS23_1} 
where each row and column is a short triangle; see \cite[prop.\ 1.4.6]{N}.
\[
\vcenter{
  \xymatrix @+0.5pc {
    a' \ar[r] \ar[d] & a_1 \ar^{ q }[r] \ar@{=}[d] & x \ar^{ i }[d] \\
    \Sigma^{ -1 }c \ar[r] \ar[d] & a_1 \ar_{ \alpha_1 }[r] \ar[d] & a_0 \ar[d] \\
    \Sigma^{ -1 }a'' \ar[r] & 0 \ar[r] & a'' \\
                    }
        }
\]
The third column and first row show that $x \in ( \Sigma^{ -1 }\cA * \cA ) \cap ( \cA * \Sigma \cA )$ whence $x \in \cA$ by Lemma \ref{lem:SMS6}(ii).  Hence the third column and first row are conflations, so $q$ is a deflation, $i$ an inflation whence $\alpha_1$ has been factored as desired.  
\end{proof}

\section{Proofs of Theorems \ref{thm:proper_tilting13} and \ref{thm:proper_tilting16}}
\label{sec:tilting}

%\begin{Definition}
%\label{def:TP_partial_order}
%The torsion pairs in a fixed abelian category $\cA$ have a partial order defined by $( \cT',\cF' ) \sqsubseteq ( \cT,\cF )$ if and only if $\cF' \subseteq \cF$.  The least and greatest elements are $( \cA,0 )$ and $( 0,\cA )$.  
%\end{Definition}

\begin{proof}
[Proof of Theorem \ref{thm:proper_tilting13}]
The values of $\Psi$ are in $\TP( \cA )$ by Lemma \ref{lem:proper_tilting11}.

The values of $\Phi$ are in $[ \cA,\Sigma \cA ]_{ n-1 }$: We must show that $\cB = \Phi( \cT,\cF ) = \Sigma \cF * \cT$ is a proper abelian subcategory of $\cC$ satisfying $E_{ n-1 }$ such that $\cA \subseteq \cB * \Sigma^{ -1 }\cB$ and $\cB \subseteq \Sigma \cA * \cA$.    All the latter conditions hold by Lemma \ref{lem:proper_tilting567}, parts (iii) and (vii), so we must show that $\cB \subseteq \cC$ is a proper abelian subcategory.  Lemma \ref{lem:proper_tilting567}, parts (iv), (v), and (vii) give that $\cB$ is an additive, extension closed subcategory of $\cC$ which satisfies condition $E_2$.  By Theorem \ref{thm:SMS7} it remains to show $\cB * \Sigma \cB \subseteq \Sigma \cB * \cB$.

We have the following string of inclusions.
\begin{align*}
  \cB * \Sigma \cB
  & \stackrel{ \rm (a) }{ \subseteq } \Sigma \cA * \cA * \Sigma^2 \cA * \Sigma \cA \\
  & \stackrel{ \rm (b) }{ \subseteq }  
  \Sigma \cA * \Sigma^2 \cA * \cA * \Sigma \cA \\
  & \stackrel{ \rm (c) }{ \subseteq }
  \Sigma^2 \cA * \Sigma \cA * \Sigma \cA * \cA \\
  & \stackrel{ \rm (d) }{ \subseteq }
  \Sigma^2 \cA * \Sigma \cA * \cA \\
  & \stackrel{ \rm (e) }{ \subseteq }
  \Sigma^2 \cB * \Sigma \cB * \Sigma \cB * \cB * \cB * \Sigma^{ -1 }\cB \\
  & \stackrel{ \rm (f) }{ \subseteq }
  \Sigma^2 \cB * \Sigma \cB * \cB * \Sigma^{ -1 }\cB 
\end{align*}
where
\begingroup
\renewcommand{\labelenumi}{(\alph{enumi})}
\begin{enumerate}
\setlength\itemsep{4pt}

  \item  is by Lemma \ref{lem:proper_tilting567}(iii),
  
  \item  is by Lemma \ref{lem:pre_proper_tilting8a}, 

  \item  is by two applications of Theorem \ref{thm:SMS7},
  
  \item  is because $\cA$ is extension closed, 

  \item  is by Lemma \ref{lem:proper_tilting567}(iii),
  
  \item  is because $\cB$ is extension closed by Lemma \ref{lem:proper_tilting567}(iv).

\end{enumerate}
\endgroup
This shows that given $c \in \cB * \Sigma \cB$, we have $c \in ( \cB * \Sigma \cB ) \cap ( \Sigma^2 \cB * \Sigma \cB * \cB * \Sigma^{ -1 }\cB )$.  Hence $c \in \Sigma \cB * \cB$ by Lemma \ref{lem:pre_proper_tilting8b}.

$\Psi\Phi = \id$:  This holds by Lemma \ref{lem:proper_tilting567}(vi).

%$n \geqslant 4$ is used in the following because we need $\cB$ to satisfy $E_3$ in order to apply Lemma \ref{lem:proper_tilting11} as described.
$\Phi\Psi = \id$:  Let $\cB \in [ \cA,\Sigma \cA ]_{ n-1 }$ be given.  Applying Lemma \ref{lem:proper_tilting11} with $\cB$ and $\Sigma \cA$ in place of $\cA$ and $\cB$ shows that $( \cB \cap \Sigma \cA,\cB \cap \cA )$ is a torsion pair in $\cB$.  Hence
\[
  \Phi\Psi( \cB )
  = \Phi( \cA \cap \cB,\cA \cap \Sigma^{ -1 }\cB )
  = \Sigma( \cA \cap \Sigma^{ -1 }\cB ) * ( \cA \cap \cB )
  = ( \cB \cap \Sigma \cA ) * ( \cB \cap \cA )
  = \cB.
\]

$\Psi$ preserves the partial orders $\sqsubseteq$ and $\preccurlyeq$: Suppose that $\cB' \preccurlyeq \cB$ in $[ \cA,\Sigma \cA ]_{ n-1 }$.  Let $\Sigma^{ -1 }b' \in \cA \cap \Sigma^{ -1 }\cB'$ be given.  Since $\cB \in [ \cA,\Sigma \cA ]_{ n-1 }$ we have $\cA \subseteq \cB * \Sigma^{ -1 }\cB$ so there is a short triangle $b_0 \xrightarrow{ \beta } \Sigma^{ -1 }b' \xrightarrow{} \Sigma^{ -1 }b_1$ with $b_i \in \cB$.  Here $\beta = 0$ by condition $E_2$ for $\cB$ since $\cB' \preccurlyeq \cB$ implies $\cB' \subseteq \cB * \Sigma^{ -1 }\cB$.  Hence $\Sigma^{ -1 }b'$ is a direct summand of $\Sigma^{ -1 }b_1$ whence $\Sigma^{ -1 }b' \in \Sigma^{ -1 }\cB$.  This shows $\Sigma^{ -1 }b' \in \cA \cap \Sigma^{ -1 }\cB$, so we have shown $\cA \cap \Sigma^{ -1 }\cB' \subseteq \cA \cap \Sigma^{ -1 }\cB$, that is, $\Psi( \cB' ) \sqsubseteq \Psi( \cB )$.  

$\Phi$ preserves the partial orders $\sqsubseteq$ and $\preccurlyeq$: Suppose that $( \cT',\cF' ) \sqsubseteq ( \cT,\cF )$ in $\TP( \cA )$, that is, $\cF' \subseteq \cF$ whence $\cT \subseteq \cT'$.  Using Lemma \ref{lem:criterion_for_torsion_pair} shows that $\cB' = \Phi( \cT',\cF' )$ and $\cB = \Phi( \cT,\cF )$ satisfy
\begin{align*}
  & \cB'
  = \Sigma \cF' * \cT'
  \subseteq \Sigma \cF * \cA 
  = \Sigma \cF * \cT * \cF
  \subseteq \Sigma \cF * \cT * \cF * \Sigma^{ -1 }\cT
  = \cB * \Sigma^{ -1 }\cB, \\
  & \cB
  = \Sigma \cF * \cT
  \subseteq \Sigma \cA * \cT'
  = \Sigma \cT' * \Sigma \cF' * \cT'
  \subseteq \Sigma^2 \cF' * \Sigma \cT' * \Sigma \cF' * \cT'
  = \Sigma \cB' * \cB',
\end{align*}
that is, $\cB' \preccurlyeq \cB$.  
\end{proof}

\begin{proof}
[Proof of Theorem \ref{thm:proper_tilting16}]
The inclusion $[ \cA,\Sigma \cA ]_n \subseteq [ \cA,\Sigma \cA ]_{ n-1 }$ is immediate from the definitions of the interval $[ \cA,\Sigma \cA ]_n$ and condition $E_n$, so we must show $[ \cA,\Sigma \cA ]_n \supseteq [ \cA,\Sigma \cA ]_{ n-1 }$ if (i) or (ii) holds.  This is clear if (i) holds so assume (ii) and let $\cB \in [ \cA,\Sigma \cA ]_{ n-1 }$ be given.  Then $\cB$ satisfies condition $E_{ n-1 }$ so to show $\cB \in [ \cA,\Sigma \cA ]_n$ it is sufficient to show $\cC( \cB,\Sigma^{ -n }\cB ) = 0$.  However, by Theorem \ref{thm:proper_tilting13} we know $\cB = \Sigma \cF * \cT$ for a torsion pair $( \cT,\cF )$ in $\cA$, so it is enough to prove the following. 
\begingroup
\renewcommand{\labelenumi}{(\alph{enumi})}
\begin{enumerate}
\setlength\itemsep{4pt}

  \item  $\cC( \Sigma \cF,\Sigma^{ -n }\Sigma \cF ) = 0$.

  \item  $\cC( \Sigma \cF,\Sigma^{ -n }\cT ) = 0$.
  
  \item  $\cC( \cT,\Sigma^{ -n }\Sigma \cF ) = 0$.

  \item  $\cC( \cT,\Sigma^{ -n } \cT ) = 0$.

\end{enumerate}
\endgroup
Since $\cF, \cT \subseteq \cA$ while $n \geqslant 4$, condition $E_n$ for $\cA$ immediately implies parts (a), (c), and (d).  For part (b) we have
\[
  \cC( \Sigma \cF,\Sigma^{ -n }\cT )
  \cong \cC( \cF,\Sigma^{ -n-1 }\cT )
  \stackrel{ \rm (e) }{ \cong } \dual \cC( \cT,\cF )
  \stackrel{ \rm (f) }{ = } 0
\]
where (e) is because $\cC$ is $(-n-1)$-Calabi--Yau, (f) is because $( \cT,\cF )$ is a torsion pair, and $\dual$ denotes $k$-linear duality.
\end{proof}

\begin{Remark}
\label{rmk:proper_tilting8}
In the situation of Theorem \ref{thm:proper_tilting13}, it follows from Lemma \ref{lem:criterion_for_torsion_pair} that $( \Sigma \cF,\cT )$ is a torsion pair in $\cB$ 
since $\cC( \Sigma \cF,\cT ) = 0$ because $\cA$ satisfies condition $E_4$, hence condition $E_1$.
\end{Remark}

\medskip
\noindent
{\bf Acknowledgement.}
We thank David Pauksztello for answering several questions.

This work was supported by a DNRF Chair from the Danish National Research Foundation (grant DNRF156), by a Research Project 2 from the Independent Research Fund Denmark (grant 1026-00050B), and by the Aarhus University Research Foundation (grant AUFF-F-2020-7-16).

\end{document}